\documentclass[a4paper,12pt]{article}

\def\DD{\displaystyle}

\usepackage{amsmath,amssymb}
\newcommand{\rmd}{\,\mathrm{d}}
\newtheorem {rema}{Remark}
\newtheorem {lemma}{Lemma}

\newtheorem {prop}{Proposition}
\newtheorem {thm}{Theorem}
\def \a   {\alpha}

\def \g   {{\gamma}}
\def \r   {{\rho}}
\def \eps {\varepsilon}
\def \dto {{\stackrel{\cal D}{\,\longrightarrow\,}}}

\def \lf   {{\lfloor}}
\def \rf   {{\rfloor}}

\def\E{{\mathbb{E\,}}}
\def\Var{{\mathbb{V}{\sf ar}\,}}
\def\P{{\mathbb{P}}}

\newcommand{\Z}{\mathbb{Z}}

\def\phi{{\varphi}}

\def\limt{\lim_{t\to\infty}}
\def\limt0{\lim_{t\to 0}}

\def\lf{\lfloor}
\def\rf{\rfloor}
\def\|{\,|\,}
\def \eps {\varepsilon}
\newcommand{\BBox}{\rule{6pt}{6pt}}
\newcommand\Cox{$\hfill \BBox$ \vskip 5mm}
\def\bn{\begin{eqnarray*}}
\def\en{\end{eqnarray*}}
\def\bnn{\begin{eqnarray}}
\def\enn{\end{eqnarray}}

\title {Forest fires on $\Z_+$ with ignition only at $0$}

\author{Stanislav Volkov\footnote{Department of Mathematics,
University of Bristol, BS8~1TW, U.K.
\newline
E-mail:~S.Volkov@bristol.ac.uk}}
\date{15 September 2009}

\begin {document}
\maketitle

\begin {abstract}
We consider a version of the forest fire model on  graph $G$,
where each vertex of a graph becomes occupied with rate one. A
fixed vertex $v_0$ is hit by lightning with the same rate, and
when this occurs, the whole cluster of occupied vertices
containing $v_0$ is burnt out. We show that when $G=\Z_{+}$, the
times between consecutive burnouts at vertex $n$, divided by $\log
n$, converge weakly as $n\to\infty$ to a random variable which
distribution is $1-\rho(x)$ where $\rho(x)$ is the Dickman
function.

We also show that on transitive graphs with a non-trivial site
percolation threshold and one infinite cluster at most, the
distributions of the time till the first burnout of {\it any}
vertex have exponential tails.

Finally, we give an elementary proof of an interesting  limit:
$\lim_{n\to\infty} \sum_{k=1}^n {n \choose k} (-1)^k \log k
-\log\log n=\g$.
\end {abstract}



\section{Introduction and results}
Consider the following forest fire model on
$\Z_+=\{0,1,2,\dots\}$. Let $\eta_x(t)\in \{0,1\}$ be the state of
site $x\in \Z_+$ at time $t\ge 0$, and we say that site $x$ is
{\it vacant} if $\eta_x=0$ and {\it occupied}, if $\eta_x=1$. The
vacant sites become occupied with rate $1$; once they are
occupied, they can only be ``burnt'' by a fire spread from a
neighbour, which reverses them to the original vacant state.
Imagine that there is a constant source of fire attached to site
$0$. Hence, whenever site $0$ becomes occupied, the whole
connected cluster of occupied sites containing $0$ is
instantaneously burnt out. Denote the process we obtain as
$\{\eta_x(t)\}$. We are interested in the dynamics of process
$\{\eta_x(t)\}$, as time passes by, under the assumption that all
sites are initially vacant, i.e. $\eta_x(0)=0$ for all $x$.

Note that our model differs from  more classical versions
presented in~\cite{VDBa} and~\cite{VDBb}, where {\it each}
occupied site can be ignited at rate $\lambda$, and then the
cluster containing this site disappears. On the other hand, our
model on $\Z_+$ turns out to be a special case of the one studied
in \cite{VDBT}, where some of the results, independently obtained
in the present paper, are also given. This covers, for example,
the recursion (\ref{equrec}), and also most of Lemma~\ref{lem_mu},
but none of the limiting statements derived in
Theorems~\ref{thm_Di} and~\ref{t3}. Forest fire models have also
been recently considered  on Erd\H{o}s-R\'enyi random graphs, see
\cite{RaTo}.

Let $T_x(i)$, $i=1,2,\dots$, be the consecutive times when site
$x$ is burnt for the $i$-th time, and let $T_x(0)=0$. Let
$\tau_n(i)=T_n(i)-T_n(i-1)$ for $i\ge 1$. We can easily show that
for a fixed $n$, $\tau_n(i)$'s are i.i.d.\ random variables; this
can be done by induction on $n$. Indeed, the times of burnouts at
$(n+1)$ depend only on $T_n$'s and the Poisson arrival process at
site $(n+1)$ itself. Since for each $j$ necessarily
$T_{n+1}(j)=T_n(i)$ for some $i$, $T_{n+1}(j)$'s are renewal
times, and hence $\tau_{n+1}(j)$'s are i.i.d.\ as well.

Now we would like to find the distribution of $\tau_{n+1}(i)$'s.
For site $0$ this is trivial as the burn-out times constitute a
Poisson process, so that
$$
 \P(\tau_0>u)= e^{-u},\ u\ge 0.
$$
Reasonably easy one can also obtain
$$
\P(\tau_1>u)=(u+1) e^{-u},
$$
so that $\tau_1$ has $\Gamma(2,1)$ distribution with density $u
e^{-u}$; similarly
$$
\P(\tau_2>u)=\frac {(2u^2+10u+7) e^{-u}+ e^{-3u}}8.
$$
From the above calculations we conclude that
 \begin{equation*}
 \begin{array}{rclrclrcl}
 \E(\tau_0)&=&1,& \E(\tau_1)&=&2, & \E(\tau_2)&=&8/3;\\
 \Var(\tau_0)&=&1, & \Var(\tau_1)&=&2,& \Var(\tau_2)&=&8/3.
 \end{array}
 \end{equation*}
Incidentally, this suggests  that  $\Var(\tau_n)=\E \tau_n$,
which, however, turns out to be incorrect, as follows from
Remark~\ref{rem_var}.

For a general $n$, let $\phi_n(t)=\E e^{t\tau_n}$ be the moment
generating function  of random variable $\tau_n$. Suppose that
sites $n$ and $n+1$ have just been burnt, and without loss of
generality reset the time to $t=0$. Let $\eta\sim \exp(1)$ be the
time till the next Poisson arrival at site $(n+1)$. The next
burnout at site $(n+1)$ will be either at time $t=\tau_n$ if
$\eta\le \tau_n$, or at a later time otherwise; in the latter case
due to the memoryless property of the Poisson process the time
between $\tau_n$ and the next burnout at $(n+1)$, denoted by
$\tilde\tau_{n+1}$, will have the same distribution as
$\tau_{n+1}$ itself. Therefore, given $\tau_n$,
 \bn
 \tau_{n+1}=\tau_n+\left\{\begin{array}{ll}
 0,&\mbox{ if }\eta\le \tau_n;\\
 \tilde\tau_{n+1},&\mbox{ if }\eta>\tau_n.\\
 \end{array}\right.
 \en
Consequently,
 \bn
 \E \left(e^{t\tau_{n+1}}\| \tau_n \right)&=&
 e^{t\tau_n} \left[(1-e^{-\tau_n})\cdot 1 + e^{-\tau_n}\cdot \E\left(e^{t\tilde\tau_{n+1}}\|\tau_n\right) \right]
 \\ &=& e^{t\tau_n}-e^{(t-1)\tau_n} + e^{(t-1)\tau_n}\cdot \E\left(e^{t\tilde\tau_{n+1}}\right)
 \en
using the fact that $\tilde\tau_{n+1}$ is independent of $\tau_n$.
Taking the expectation on both sides, we obtain
 \bn
 \phi_{n+1}(t)=\phi_{n}(t)-\phi_{n}(t-1)+\phi_{n}(t-1)\phi_{n+1}(t)
 \en
whence
 \bnn\label{eqphirec}
 \phi_{n+1}(t)=\frac{\phi_{n}(t)-\phi_{n}(t-1)}{1-\phi_{n}(t-1)}.
 \enn
Let
$
 u_n(t)=\phi_n(t)-1.
$
Then
 \bnn\label{equrec}
 u_{n+1}(t)=-\frac{u_{n}(t)}{u_n(t-1)}.
 \enn
Since $\phi_0(t)=\int_0^{\infty} e^{tx} e^{-x}\rmd x =1/(1-t)$
yielding $u_0(t)=t/(1-t)$, we can easily iteratively compute
$u_n(t)$. For example,
$$
u_1(t)=\frac{t(2-t)}{(1-t)^2},\
u_2(t)=\frac{t(2-t)^3}{(1-t)^3(3-t)}
$$
which is consistent with our previous calculations of the
distributions of $\tau_1$ and $\tau_2$.
\begin{lemma}\label{lem_u}
For $n=1,2,...$
 \bn
 u_{n-1}(t)=\frac{t\cdot (2-t)^{{n\choose 2}} \cdot  (4-t)^{{n\choose 4}}  \cdot (6-t)^{{n\choose 6}} \dots }
  {(1-t)^{{n\choose 1}} \cdot  (3-t)^{{n\choose 3}} \cdot  (5-t)^{{n\choose 5}}  \dots}
 \en
with the convention that ${n\choose k}=0$ whenever $k>n$. Thus
$\tau_n$ is a mixture of Gamma random variables with the moment
generating function
 \bn
 \phi_{n-1}(t)&=&1+t\prod_{k=1}^{n} (k-t)^{(-1)^k{n\choose {k}} }
  \en
defined for all $t<1$.
\end{lemma}
{\sl Proof.} By induction, using  (\ref{equrec}) and the fact that
${n\choose k}+{n\choose {k-1}}={{n+1}\choose k}$. \Cox

\begin{lemma}\label{lem_mu}
Let $\mu_n=\E \tau_n$. Then for $n=1,2,...$
 \bn
  \log \mu_{n-1}&=&\sum_{i=1}^n {n\choose i} (-1)^{i} \log i;\\
  \\
  \E \left(\tau_{n-1}^2\right)&=&2 \mu_{n-1} \sum_{i=1}^n {n\choose i}
  \frac{(-1)^{i}}{i}.
 \en
 Moreover,
 \bn
  \lim_{n\to\infty} \log\left(\frac{\mu_n}{\log n}\right)=\g
 \en
where $\g=0.577...$ is the Euler constant.
\end{lemma}

The following two lemmas will be proved in Section~\ref{sec_eg}.

\begin{lemma}\label{lem_com1}
Let
 \bn
 A_{n-1}&=&\sum_{i=1}^n {n\choose i} (-1)^{i} \log
 i
 \en
then $\lim_{n\to\infty} (A_n-\log\log n)= \g $, where
$\g=0.577...$ is the Euler constant.
\end{lemma}
\begin{rema}
After this paper has been written, we learned (B\'alint T\'oth,
personal communications) that the above limit is in fact derived
in~\cite{FLAJ}, Theorem 4, though no explicit proof was given
there. Thus we shall give a reasonably short and {\sf elementary}
self-contained proof of this convergence.
\end{rema}

\begin{lemma}\label{lem_com2}
Let
 \bn
 a(n,m)=\sum_{k=1}^n {n \choose k}  \frac{(-1)^{k+1}}{k^m}.
 \en
Then
 \bn
  & (a) & a(n,m)=  \sum_{1\le i_1\le i_2\le\dots\le i_m\le n} \frac 1{i_1 i_2\dots
  i_m};
  \\ \\
  & (b) & a(n,m)\le (\log n +1)^m \mbox{ for all } n\ge 1;
  \\ \\
  & (c) & a(n,m)=\frac{\log^m n}{m!}+O(\log^{m-1} n) \mbox{  for a fixed $m$ as
  }n\to\infty.
 \en
\end{lemma}
\begin{rema}
The quantities $a(n,m)$ are closely related to the Stirling
numbers of the second kind:
$$
\left\{\begin{array}{c}n\\
m\end{array}\right\}=\frac 1{n!}\sum_{k=0}^n {n\choose
k}(-1)^{n-k} k^m,
$$
and up to a coefficient of proportionality coincide with
``negative-positive'' Stirling numbers in \cite{BRAN}, see
equations (68) and (78) there.
\end{rema}

{\sl Proof of Lemma~\ref{lem_mu}.} The first part follows
immediately from Lemma~\ref{lem_u} and the properties of
moment-generating functions; the second part follows from
Lemma~\ref{lem_com1}. \Cox

\begin{rema}\label{rem_var}
Lemma~\ref{lem_mu} together with part (c) of Lemma~\ref{lem_com2}
yield that
 \bn
  \lim_{n\to\infty} \E\left(\frac{\tau_n}{\log n}\right)&=&e^{\g}\approx 1.78\dots \\
  \lim_{n\to\infty} \E\left(\frac{\tau_n}{\log n}\right)^2 &=&2e^{\g}\approx
  3.56\dots
 \en
whence for large $n$, $\Var(\tau_n)\ne \E \tau_n$.
\end{rema}

\begin{thm}\label{thm_conv}
Let $\xi_n=\tau_n/\log n$. Then as $n\to\infty$
 \bn
  \xi_n\dto \xi
 \en
(meaning convergence in distribution) where $\xi$ is a random
variable with mean $\E \xi =\g'$, and the moment generating
function
\begin{eqnarray}\label{eqlimchf}
  \phi_{\xi}(s) &\equiv & \E e^{s\xi }=1+\exp\left\{{\rm Ei}(s)\right\}.
\end{eqnarray}
Here $\g'=e^{\g}=1.781...$ and
 \bn
 {\rm Ei}(s)=\int_{-\infty}^s \frac{e^x}{x} \rmd x=\g+\log s +\sum_{m=1}^{\infty}
\frac{s^m}{m\cdot m!}
 \en
is the exponential integral (understood in terms of the Cauchy
principal value; see~\cite{AS}, Section 5.1 and formula 5.1.10).
\end{thm}
{\sl Proof.} Observe that
 \bn
\log \left[\frac{u_{n-1}(t)}{t \mu_{n-1}}\right]&=& \sum_{k=1}^n
{n \choose k} (-1)^{k}\log\left(1-\frac tk\right)=
\sum_{m=1}^{\infty} \frac{t^m}{m} \left[ \sum_{k=1}^n {n \choose
k} \frac{(-1)^{k-1}}{k^m}\right]
 \\  \\ &=& \sum_{m=1}^{\infty}  \frac{a(n,m)\,t^m}{m}
 \en
where
 \bn
 a(n,m)=\sum_{k=1}^n {n \choose k}
 \frac{(-1)^{k+1}}{k^m}.
 \en
For the moment generating function of $\xi_{n-1}$ we have
 \bn
\log(\E e^{s\xi_{n-1}} -1)=\log u_{n-1}(s/\log n)&=&\log s +\log
\frac{\mu_{n-1}}{\log n} +
 \sum_{m=1}^{\infty} \frac{a(n,m)}{\log^m n} \frac{s^m}{m},
 \en
consequently for any $N\ge 1$, using part (b) of
Lemma~\ref{lem_com2},
 \bn
 \Delta_{n-1}(s)&:=&\left|\log(\E e^{s\xi_{n-1}} -1)-{\rm Ei}(s)\right| \\
 &=& \left|\log(\E e^{s\xi_{n-1}} -1)-\g-\log s -\sum_{m=1}^{\infty} \frac{s^m}{m\cdot m!}\right| \\ \\
 &\le& \left|\log \frac{\mu_{n-1}}{\log n}-\g\right| +
 \sum_{m=1}^{N} \left|\frac{a(n,m)}{\log^m n}-\frac 1{m!}\right| \frac{s^m}{m}
 \\ \\
 &+& \sum_{m=N+1}^{\infty}  \frac{s^m}{m\cdot m!}
 + \sum_{m=N+1}^{\infty} \left(1+\frac 1{\log n}\right)^m
\frac{s^m}{m}.
 \en
Fix an $\eps>0$. Assuming $|s|\le 1/2$, we can choose $N$ so large
that the last two summands are smaller than $\eps/2$ each. Now for
a {\it fixed} $N$ by Lemma~\ref{lem_mu} and  part (c) of
Lemma~\ref{lem_com2} the first two terms of the RHS of
$\Delta_{n-1}(s)$ go to $0$ as $n\to\infty$. Consequently,
$\limsup_{n\to\infty} \Delta_n(s)\le \eps$. Since $\eps$ is
arbitrary, we conclude that for $|s|\le 1/2$
$$
\lim_{n\to\infty} \E e^{s\xi_{n-1}}=1+ \exp\left\{{\rm
Ei}(s)\right\}.
$$
By Theorem 3 in~\cite{CUR}, if the sequence  of moment-generating
functions corresponding to random variables $\xi_n$ converges
point-wise to a limit function $\phi_{\xi}(s)$ on {\it some}
interval around $0$, then there is a random variable $\xi$ such
that $\xi_n \to \xi$ in distribution and $\phi_{\xi}(s)$ is its
moment generating function. This finishes the proof. \Cox

\begin{thm}\label{thm_Di}
Random variable $\xi$ defined in Theorem~\ref{thm_conv} has the
density function $f(x)$ and the survival function $\r(x)=\P(\xi>
x)$ satisfying
 \begin{equation}\label{eq_eq4f}
 \begin{array}{rcll}
 f(x)&=&0, & x\le 1; \\
 \frac{\rmd }{\rmd x} (x f(x))&=& - f(x-1), & x> 1,
 \end{array}
 \end{equation}
and
 \begin{equation}\label{eq_eq4F}
 \begin{array}{rcll}
 \r(x)&=&1, & x\le 1;\\
  x\r'(x)&=&-\r(x-1), &x>1,
 \end{array}
 \end{equation}
so that $\r(x)$ is the Dickman function.
\end{thm}
{\sl Proof.} Let us denote by $\psi(t)=\E e^{it
\xi}=\phi_{\xi}(it)$, then we have
$$
t\psi'(t)=\psi(t)e^{it}-e^{it}.
$$
Using {\it formally} the inversion formula and the fact that for a
random variable $Y\equiv 1$, $\E e^{itY}=e^{it}$, we have
 \bn
\frac{1}{2\pi} \int (t\psi(t))'e^{-itx} \rmd t&=&\frac{1}{2\pi}
\int [\psi(t)e^{it}+\psi(t)] e^{-itx} \rmd t-\delta_{x-1}\\ \\
 &=&\frac{1}{2\pi} \int \psi(t) e^{-it(x-1)} \rmd t+\frac{1}{2\pi}
\int \psi(t) e^{-itx} \rmd t-\delta_{x-1}
 \en
(where $\delta_x$ denotes the Dirac delta-function.) Using
integration by parts on the left (again, formally) we have
 \bn
ix \frac{1}{2\pi} \int t \psi(t)e^{-itx} \rmd t
 &=&-x\frac{\rmd}{\rmd x} \left[ \frac{1}{2\pi} \int \psi(t) e^{-itx} \rmd t\right]=-x\frac{\rmd}{\rmd x} f(x)
 \\ \\ &=& f(x-1)+f(x)-\delta_{x-1}
 \en
yielding $(x f(x))'=\delta_{x-1}-f(x-1)$. Integrating this
equality from $-\infty$ to $x$, and denoting $F(x)=\P(\xi\le x)$,
we obtain
$$
xF'(x)=1_{x\ge 1}-F(x-1)
$$
implying $\r(x)=1-F(x)$ satisfies $x\r'(x)=-\r(x-1)$ for $x\ge 1$
as required.

To prove the above results rigorously, first observe that the
Dickman function $\r(u)$ has the following properties: (1) it is
positive and decreasing on the $[1,\infty)$; (2) it is infinitely
differentiable on $[0,\infty]$ except at integer points; (3)
$\r(u)\le 1/\Gamma(u+1)$ for $u\ge 1$ (see e.g.\  \cite{Tenen} for
its properties). Consequently,
$$
F(u):=\left\{\begin{array}{ll}
 0, &u<1;\\
 1-\rho(u), &u\ge 1
\end{array}\right.
$$
is the cumulative distribution function of some continuous random
variable $\zeta$ which density is supported on $[1,\infty)$.
Multiplying the second equation in (\ref{eq_eq4F}) by $te^{tx}$
and integrating, we obtain
 \bnn\label{eq_LTD}
\int_{1}^{\infty}  tx F'(x) e^{tx}\rmd x= \int_{1}^{\infty} (1-
F(x-1))te^{tx}\rmd x
 \enn
Integrating by parts the RHS of (\ref{eq_LTD}), we have
 \bn
 t\lim_{x\to\infty} \rho(x-1) e^{tx} - [1-F(0)]e^t+ \int_{1}^{\infty} e^{tx} F'(x-1) \rmd x
 \\ \\
 =0-e^t+e^t \int_{0}^{\infty} e^{tx} F'(y) \rmd y =e^t  (\phi_{\zeta}(t)-1)
 \en
where $\phi_{\zeta}(t)=\E e^{t\zeta}$ is the moment generating
function of $\zeta$. On the other hand, the LHS of (\ref{eq_LTD})
equals
 \bn
t\, \frac{\rmd}{\rmd t}\int_{1}^{\infty}   F'(x) e^{tx}\rmd
x=t\phi_{\zeta}'(t).
 \en
This yields $t\phi_{\zeta}'(t)=e^t (\phi_{\zeta}(t)-1)$ and
$\phi_{\zeta}(0)=1$, a general solution to which has a form
 \bn
  \phi_{\zeta}(t)=1+C_1 s  \exp\left(\sum_{m=1}^{\infty} \frac {t^m}{m\cdot m! } \right)
 \en
for some constant $C_1$. To identify $C_1$, we will use the fact
that $\phi_{\zeta}(-z)=\E e^{-z\zeta}\to 0$ as $z\to \infty$ (as
$\zeta\ge 0$).  Using Taylor expansion for $e^{t}$ we obtain
 \bn
  \phi_{\zeta}(-z)=1-C_1 z  \exp\left( - \int_0^{z} \frac{1-e^{-t}}{t}\rmd t   \right).
 \en
Now, formulas 5.1.1 and 5.1.39 in~\cite{AS} for the function
$E_1(z)$ give
 \begin{equation}\label{eqE1}
 \DD
 \begin{array}{rcl}
 E_1(z)&=&\int\limits_{z}^{\infty}\frac{\DD e^{-t}}{\DD t} \rmd t,\\
 \int\limits_{0}^{z}\frac{\DD 1-e^{-t}}{\DD t} \rmd t&=&E_1(z)+\log z +\g
 \end{array}
 \end{equation}
yielding
 \bn
  \phi_{\zeta}(-z)=1-C_1   \exp\left(-\gamma - \int_z^{\infty} \frac{e^{-t}}{t}\rmd t   \right).
 \en
Since the integral goes to $0$ as $z\to\infty$, we conclude that
$C_1=e^{\gamma}$. Thus $\phi_{\zeta}$ coincides with the
expression given by (\ref{eqlimchf}) and by the uniqueness
theorem, $\xi$ must have the same distribution as $\zeta$, from
which the Theorem follows.
 \Cox

Here are a few observations about the distribution of $\xi$.
Trivially we have $F(x)=0$ for $x\le 0$; thus using
(\ref{eq_eq4F}) for $0\le x\le 1$ we have $F'(x)=0$ whence
$$
F(x)=0, \\ 0\le x\le 1
$$
as well. Consequently, for $1\le x\le 2$, we have $x F'(x)=1$ so
that
$$
F(x)=\log x, \\ 1\le x\le 2.
$$
Therefore, by induction we can obtain piece-wise smooth density
function of $\xi$:
 \bn
 f(x)=\left\{\begin{array}{ll}
 0, &x\le 1;\\
 1/x, &1<x\le 2;\\
 \frac{1-\log(x-1)}{x}, &2<x\le 3;\\
...
 \end{array}\right.
 \en
Unfortunately, there is no explicit formula in elementary
functions for $f(x)$ on an interval $[n,n+1]$ for $n\ge 2$.

Our next statement deals with residual waiting times for the
renewal process generated by consecutive burnouts at site $n$.
\begin{prop}
Let $\eta_{t,n}$ be the time till the next burnout at site $n$
after time $t>0$. Then $\eta_{t,n}/\log n\dto \bar\eta_n$ as
$t\to\infty$ and $\bar\eta_n\dto \bar\eta$ as $n\to\infty$, where
$\bar\eta$ has a generalized Dickman distribution GD(1), see
\cite{PeWa}, i.e.\ the same distribution as $U_1+U_1 U_2+U_1 U_2
U_3+\dots$ with $U_i$ being i.i.d.\ uniform $[0,1]$ random
variables.
\end{prop}
{\sl Proof.} As we already know, the times between consecutive
burnouts $\tau_n^{(i)}$, $i=1,2,\dots$ are i.i.d.\ and  have a
common distribution of $\tau_n$. Let
$\xi_n^{(i)}=\tau_n^{(i)}/\log n$, and let $F_n(\cdot)$ be the
common cumulative distribution function of $\xi_n^{(i)}$, which is
the same as for the random variable $\xi_n$ defined in
Theorem~\ref{thm_conv}. As it is well-known, see e.g.~\cite{DUR},
Chapter 3.4, the residual waiting times for the renewal process
generated by $\xi_n^{(i)}$  converge in distribution to a
non-negative random variable $\bar\eta_n$ such that
 \bn
 \P(\bar\eta_n\le x)=\frac{1}{\E \xi_n} \int_{0}^{x} (1-F_n(u))\rmd
 u, \ \mbox{ for all }x\ge 0.
 \en
(We need to verify that the distribution $F_n$ is non-arithmetic,
however this easily follows from the fact that $\tau_n$  is a
continuous random variable, which is a mixture of Gamma
distributions, as implied by Lemma~\ref{lem_u}.)

Let $F$ be the cumulative distribution function of $\xi$, as
defined in the proof of Theorem~\ref{thm_Di}, and $\bar\eta$ be a
non-negative random variable such that
 \bnn\label{eq_GD1}
 \P(\bar\eta\le x)=\frac{1}{\E \xi} \int_{0}^{x} (1-F(u))\rmd
 u =e^{-\gamma}  \int_{0}^{x} \rho(u) \rmd
 u , \ \mbox{ for all }x\ge 0.
 \enn
Then, since $\int_{0}^{x} (1-F_n(u))\rmd u\le \int_{0}^{\infty}
(1-F_n(u))\rmd u=\E \xi_n$,
 \bn
 |\P(\bar\eta_n\le x)-\P(\bar\eta\le x)|
 &\le& \left|\frac{1}{\E \xi} \int_{0}^{x} (F_n(u)-F(u))\rmd u\right|
 \\ & \ {}\ {}\ & + \left|\frac{1}{\E \xi}-\frac{1}{\E \xi_n}\right| \int_{0}^{x} (1-F_n(u))\rmd u
 \\
 &\le& e^{-\gamma}\int_{0}^{x} \left|F_n(u)-F(u)\right|\rmd u+\left|\frac{\E\xi_n}{\E
 \xi}-1\right|\to 0
 \en
where the first summand tends to $0$ by the dominated convergence
theorem since $F_n(x)\to F(x)$ pointwise by
Theorems~\ref{thm_conv} and~\ref{thm_Di}, and the second one
vanishes because of Lemma~\ref{lem_mu}. Therefore
$\bar\eta_n\dto\bar\eta$.

We finish the proof by noting that the distribution in
(\ref{eq_GD1}) coincides with the distribution of
$\sum_{i=1}^{\infty} \prod_{j=1}^i U_j$, see \cite{JMFC}.
 \Cox

We conclude by noting that similar distributions (called {\it
Dickman-type} distributions) show up in some other probabilistic
models, including e.g.\ minimal directed spanning trees  as well
as number-theory related problems, see \cite{PeWa} and references
therein. Another interesting application is in economics, related
to plot-size distributions: see~\cite{Exner}, formula (4), which
is identical to that for $\bar\eta$.

\section{Generalizations}
One can consider a similar forest fire model on an arbitrary
connected locally-finite graph $G$ with the vertex set $V(G)$ and
one special vertex $v_0\in V(G)$ which is called {\it the origin}.
Let $\eta_x(t)\in \{0,1\}$ be the state of site $x\in V(G)$ at
time $t\ge 0$; again the site $x$ is vacant (occupied resp.) if
$\eta_x=0$ ($\eta_x=1$ resp.). Vacant sites become occupied at
rate $1$; they remain occupied until they are burnt out, which
makes them vacant again. For definiteness, at time $0$ all sites
are vacant. As before, only site $v_0$ is constantly hit by
lightning, hence whenever it becomes occupied all the sites in the
cluster of occupied sites containing $v_0$ are instantaneously
burnt out.

Unfortunately, this model turns out to be not so interesting,
provided that the critical percolation threshold $p_c$ for site
percolation on $G$ is strictly smaller than $1$, which is true on
many graphs. Recall that if $\theta_{v_0}(p)=\theta(p)$ denotes
the probability that site $v_0$ belongs to an infinite cluster of
occupied sites given that {\it each} site is independently
occupied with probability $p$, then the critical percolation
threshold is defined by
$$
p_c=\sup\{p:\ \theta(p)=0\}
$$
(see for example \cite{GRIM}).

We claim that if $p_c<1$, then in our forest fire model infinitely
many sites can be burnt in a finite time. Indeed, fix a
$p\in(p_c,1)$, and let
 \bnn\label{eq_Sofp}
S=S(p)=-\log(1-p).
 \enn
Then with probability at least $ \frac{1-p}{p} \, \theta(p)>0 $
site $v_0$ becomes occupied in time exceeding $S$ (at which point
it is immediately burnt), and by that time there will be already
an infinite cluster attached to $v_0$, so that it will burn some
arbitrarily far away vertices.

As it is well-known, on many graphs ($\Z^d$, $d\ge 2$, regular
trees, some others) the number of infinite occupied clusters can
be either $0$, $1$, or $\infty$ (see \cite{NEWM}); also it is
known that on $\Z^d$, $d\ge 2$, and some infinite Cayley graphs
(but not a regular tree) the infinite cluster, whenever present,
must be unique; see~\cite{BuKe}, \cite{HaPe}, and also Chapter~8.9
in~\cite{GRIM} and Theorem~4 in Chapter 5.1 in~\cite{BoloRio}.
Additionally, suppose that the graph is {\it transitive}, that is
to say that graph $G$ viewed from any vertex $v\in V(G)$  is
isomorphic to graph $G$ viewed from $v_0$; this in turn would
imply  using the FKG inequality for the connectivity function
(\cite{GRIM}, Chapter 8.5) that the probability that an arbitrary
chosen vertex $v$ is burnt out in time $S$ exceeds $(1-\g)(1-p)$
where
 \bnn\label{eq_burn}
 \gamma:=1-\frac{\theta(p)^2 }{p} \in(0,1).
 \enn
We can generalize this argument as follows.
\begin{thm}\label{t3}
Suppose that  graph $G$ is connected, transitive, the critical
point for the site percolation $p_c=p_c(G)<1$ and that there can
be at most one infinite cluster on $G$. Fix an arbitrary $v\in
V(G)$ and let $\eta$ be the time till its first burnout in our
forest fire model. Then for any $p\in (p_c,1)$ and $j=0,1,2,\dots$
 \bnn\label{eq_eta}
 \P\left(\eta > x\right)\le \g^{-1} \left[x(1-p)+1\right]  e^{  -\lambda x}
\  \text{ for all $x>0$}
 \enn
 where  $\g$ is given by (\ref{eq_burn}), and $\lambda=\lambda(\gamma)>0$ is the
smallest positive solution of
 \bnn\label{eq_lamb}
 \phi(\lambda)=\gamma^{-1}
 \enn
with
$ \phi(t)=\left[1-\frac{t}{(1-p)^{1-t}}\right]^{-1}=  [1 -t
e^{S(1-t)}]^{-1}$
and $S$ being defined by (\ref{eq_Sofp}).
\end{thm}
\begin{rema}
The function $\phi(t)$  satisfies the following properties:
\begin{itemize}
\item $\phi(0)=1$;
\item $\phi(t)$ is positive and finite on $[0,t_{\max})$
where $t_{\max}=t_{\max}(S)$ is the smallest positive solution of
$1=t e^{S(1-t)}$, that is
$$
t_{\max}(S)=\left\{\begin{array}{ll}
 1, &\text{for }S\le 1,\\
-{\sf LambertW}(-Se^{-S})/S &\text{for }S>1 \end{array}\right.
$$
where ${\sf LambertW}$ is the Lambert W function;
\item $t_{\max}\le S^{-1}$ and hence $\phi'(t)\propto
(1-tS)>0$ for $t<t_{\max}$ (easy to check);
\item $\phi(t)\uparrow +\infty$ as $t\uparrow t_{\max}$.
\end{itemize}
Therefore, the solution to (\ref{eq_lamb}) indeed exists for any
$0<\g<1$.
\end{rema}

 \noindent
 {\sl Proof of Theorem~\ref{t3}.} As we have already established, the
probability that an arbitrary vertex $v$ is burnt out in time $S$
is at least (\ref{eq_burn}); this would be obviously also true
even if some of the vertices $v\in V(G)\setminus\{v_0\}$ were
already occupied at time $0$. Denote by $T(1), T(2),\dots$ the
times of ignitions of vertex $v_0$, set $T(0)=0$ and let
$\tau(n)=T(n)-T(n-1)$ be the (exponentially(1) distributed) times
between consecutive burnouts. Let $N=N(x)$ be the number of
intervals $\tau(i)$ of length at least $S$ entirely lying inside
$[0,x]$, that is
$$
 N(x)={\sf card}\{i:\ \tau_i\ge S,\ T_i\le x\}
$$
To get a handle on $N(x)$, we will use the renewal theory
approach. Let
 \bn
  i_0&=&0,
  \\
  i_k&=&\min\{i> i_{k-1}:\ \tau_i\ge S,\ \tau_j<S\ \forall
  j\in (i_{k-1},i)\},\ k=1,2,3,\dots
 \en
Then $T(i_k)$ form a renewal process, and
$$
 N(x)=\max\{k:\ T(i_k)\le x\}=\max\{k:\ \nu_1+\nu_2+\dots+\nu_k\le x\}
$$
where  $\nu_k:=T(i_{k+1})-T(i_k)$ are i.i.d.\ random variables,
and if $\phi_\nu(t)=\phi(t)$ denotes its moment generating
function which we will need later, then, by conditioning on
$\tau_1$ and using the memoryless property, we obtain
 \bn
  \phi(t)&=&\E e^{tT(i_1)}=
  \E \left[e^{tT(i_1)} 1_{\tau_1\le  S}\right]
  +\E \left[e^{tT(i_1)} 1_{\tau_1>  S}\right]\\
  &=& \E \int_0^{S} e^{tu+T(i_1)}e^{-u} \rmd u+\int_S^{\infty} e^{tu}e^{-u} \rmd u
  \\
  &=& \frac
  1{1-t}\left[\left(1-e^{-(1-t)S}\right)\phi(t)+e^{-(1-t)S})\right]
  \en
yielding
 \bn
 \phi(t)=\frac 1{1-t e^{S(1-t)}}
 \en
which is defined for all $t<t_{\max}$. In particular, $\E
\nu=\phi'(0)=e^S$, and thus we expect $N(x)$ to be typically
around $x e^{-S}=x(1-p)$.

On the other hand, by the arguments preceding the statement of the
Theorem, conditioned on $N(x)$, the probability that $v$ has {\it
not} been burnt out in time $x$ is smaller than $\g^{N(x)}$, hence
 \bn
 \P(\eta>x)\le \E \g^{N(x)}=\sum_{n=0}^{\infty} \g^n \P(N(x)=n).
 \en
We split the sum above into two parts and estimate it as follows:
 \bnn\label{eq_sum_est}
 \sum_{n=0}^{\infty} \g^n \P(N(x)=n)
 &\le&
   \sum_{n=0}^{\lf x(1-p)\rf -1} \g^n \P(N(x)=n)
 + \sum_{n=\lf x(1-p)\rf}^{\infty} \g^n \P(N(x)=n)  \nonumber \\
 &\le&
 \sum_{n=0}^{\lf x(1-p)\rf  -1} \g^n \P(N(x)\le n)
                            + \g^{\lf x(1-p)\rf}\P\left(N(x)\ge \lf x(1-p) \rf \right)  \nonumber \\
 &\le &
 \sum_{n=0}^{\lf x(1-p)\rf -1} \g^{n} \P(\nu_1+\dots+\nu_{n+1}\ge
 x)+\g^{x(1-p)-1}  \\
 &\le &
  \g^{-1} \left[x(1-p)   \max_{m\in \{1,\dots, \lf x{(1-p)}\rf \}} \g^m \P(\nu_1+\dots+\nu_{m}\ge x)
 +\g^{x(1-p)}\right]  \nonumber .
 \enn
From Markov inequality, we have for any $t>0$
 \bn
 \g^{m} \P(\nu_1+\dots+\nu_{m}\ge x) \le e^{\Lambda(t,m)} \text{ where }\Lambda(t,m)=m \log \g +m \log \phi(t) -t x.
 \en
We will bound $\log \left[\g^{m} \P(\nu_1+\dots+\nu_{m}\ge
x)\right]$ by $\max_{0\le m\le x(1-p)} \min_{t>0} \Lambda(t,m)$.
From well-known properties of the MGF we know that $\log \phi(t)$
and hence $\Lambda(t,m)$ is convex in $t$, therefore the latter
achieves a unique minimum at point $t^*=t^*(x/m)$ where
$t^*(\alpha)$ solves the equation
$$
\frac{\phi'(t^*(\a))}{\phi(t^*(\a))}=\a.
$$
Also, for $m\le  x(1-p)$ we have $t^*\ge  0$ as
$\partial\Lambda(t,m)/\partial
t\|_{t=0}=m\phi'(0)/\phi(0)-x=m(1-p)^{-1}-x\le 0$, yielding
$\min_{t\ge 0} \Lambda(t,m)=\Lambda(t^*(x/m),m)$. Additionally,
$t^*\left(\frac 1{1-p}\right)=0$, $t^*(\a)$ is increasing in $\a$
as $\rmd \log \phi(t)/\rmd t$ is increasing, and it is easy to
check in our case $\phi(t^*(\a))\to\infty$ as $\a\to\infty$.

On the other hand,
 \bn
 \frac{\rmd \Lambda(t^*(x/m),m)}{\rmd
 m}=\log\g+\log\phi(t^*)+\left[m\frac{\phi'(t^*)}{\phi(t^*)}-x\right]\frac{\rmd t^*(x/m)}{\rmd
 m}=\log\left[\g\phi\left(t^*\left(\frac xm\right)\right)\right].
 \en
The RHS of this expression decays in $m$; moreover as $m\downarrow
0$, $\frac xm \to +\infty$ resulting in $\phi(t^*(x/m))\to
+\infty$ and $\left. \frac{\rmd \Lambda(t^*(x/m),m)}{\rmd
m}\right|_{m=0}=+\infty$. At the same time, for $m=x(1-p)$ we have
$t^*(x/m)=0$ hence $\left. \frac{\rmd \Lambda(t^*(x/m),m)}{\rmd
m}\right|_{m=x(1-p)}=\log \g <0$. Therefore, the maximum of
$\Lambda(t^*(x/m),m)$ is achieved at some intermediate $m$ and
this maximum equals $-\lambda x$ where $\lambda=t^*(x/m)$ solves
$\frac{\rmd \Lambda(t^*(x/m),m)}{\rmd m}=0$, i.e.\
equation~(\ref{eq_lamb}). Finally, observe that
 \bn
 \g^{x(1-p)}&=&\exp\{\Lambda(0,x(1-p))\}= \exp\{\Lambda(t^*((1-p)^{-1}),x(1-p))\}
 \\
 &\le&     \exp\left\{ \max_{0\le m\le x(1-p)} \Lambda(t^*(x/m),m)  \right\}
 = e^{-\lambda x}.
 \en
Now (\ref{eq_sum_est}) yields (\ref{eq_eta}).
 \Cox

\section{Proofs of the combinatorial results}\label{sec_eg}

{\sl Proof of Lemma~\ref{lem_com1}.} Observe that
 \bn
 A_{n-1}&=&\sum_{i=2}^n {n\choose i} (-1)^{i} \log i=
 \sum_{i=2}^n {n\choose i} (-1)^{i} \left[\log \frac 21+\log \frac 32 +\dots+\log
 \frac{i}{i-1}\right]
 \\ \\
 &=&\sum_{i=2}^n  (-1)^{i} \left[{n\choose i}-{n\choose {i+1}}+{n\choose {i+2}}-\dots\pm {n\choose n}  \right]\log \frac{i}{i-1}
 \\ \\
 &=&\sum_{i=2}^n  (-1)^{i} {n-1 \choose i-1}\log \frac{i}{i-1}
  = \sum_{k=1}^{n-1}  (-1)^{k-1} {n-1 \choose k}\log
  \frac{k+1}{k},
 \en
hence
 \bn
 A_{n}=\sum_{k=1}^{n}  (-1)^{k-1} {n \choose k}\log \frac{k+1}{k}.
 \en
To estimate the above quantity, we use the partial fractions
method the way it is employed in \cite{SONa}, equation (8), and in
\cite{SONb}, Example 5.8,
  \bn
  \frac{n!}{x(x+1)\dots(x+n)}=\frac 1x-\sum_{k=1}^n {n\choose k}(-1)^{k-1} \frac
  1{x+k}.
  \en
Consequently,
 \bnn\label{eqintrep}
  A_n= \int_0^1 \left[\frac 1x - \frac{n!}{x(x+1)\dots(x+n)} \right]\rmd x.
 \enn
(In fact, there is yet another formula for $A_n$ in~\cite{PRUD},
5.5.1, saying that
 \bn
   \sum_{k=1}^{n}  (-1)^{k} {n \choose k}\log
   \frac{k+a}{k+b}=-\log \frac ab
   +\int_0^1(t^{a-1}-t^{b-1})(1-t)^n \frac{\rmd t}{\log t}
 \en
hence $\log A_n=\lim_{b\downarrow 0} \left[\log b^{-1}
   -\int_0^1(1-t^{b-1})(1-t)^n (\log t)^{-1} \rmd t \right]$.
Unfortunately, we could not estimate this limit and hence decided
to work directly with (\ref{eqintrep}).)

Let us rearrange (\ref{eqintrep}) as follows:
 \bn
 A_n= \int_0^1 \left[\frac 1x - \frac{1}{x(1+x/1)(1+x/2)\dots(1+x/n)} \right]\rmd x
 \en
Using standard Taylor series expansion for $|x|<1$ we have
 \bn
 \log \left(\left(1+x\right)\left(1+\frac x2\right)\dots\left(1+\frac xn\right) \right)
 = \sum_{m=1}^{\infty} \frac{x^m}{m} (-1)^{m-1} H_{n,m}
 \en
where $H_{n,m}=\sum_{k=1}^n k^{-m}$ are the generalized harmonic
numbers.
Moreover
 \bnn\label{eq_Harmn}
  H_n\equiv H_{n,1}&=&\g+\log n + \frac 1{2n}+O(n^{-2}),\\
 H_{n,m}&=&\zeta(m)-\frac 1{(m-1)n^{m-1}}+O(n^{-m}), \
 m=2,3,\dots\nonumber
 \enn
with $\zeta(s)=\sum_{k=1}^{\infty} k^{-s}$ being the Riemann zeta
function and $\g=0.577...$ the Euler constant. The first equality
in~(\ref{eq_Harmn}) follows from the asymptotic for the Digamma
function $\psi(x)=\frac{\rmd \log \Gamma(x)}{\rmd x}$ (see 6.3.2
and 6.3.18 in~\cite{AS}) while the second one is an elementary
consequence of the fact that
 $$
 \zeta(m)-H_{n,m}=\sum_{k=n+1}^{\infty} \frac 1{k^m}, \
 \text{while}\
\int_{n+1}^{\infty} \frac{\rmd x}{x^m}<\sum_{k=n+1}^{\infty} \frac
1{k^m}<\int_{n}^{\infty} \frac{\rmd x}{x^m}.
 $$

By changing the variables $x=y/\log n$ in the integral, we obtain
 \bn
 A_n= \int_0^{\log n} \frac{\rmd y}{y}\left[1 - \frac{1}{\exp(B_n(y))} \right]
 \en
where
 \bn
 B_n(y)= y+\frac{y}{\log n}\left(\g+\frac 1{2n}+O\left(n^{-2}\right)\right)-\frac{y^2}{2(\log n)^2}\left(\zeta(2)-\frac 1n
 +O(n^{-2})\right)+\dots
 \en
Integrating separately on $[0,1]$ and $[1,\log n]$ we obtain
 \begin{eqnarray*}\begin{array}{rclrlrlrl}
 A_n&=&
  \int_1^{\log n} \frac{\rmd y}{y}
  &-& \int_1^{\log n} \frac{e^{-B_n(y)}\rmd y}{y}
  &+& \int_0^{1} \frac{\rmd y}{y}\left[1 - e^{-B_n(y)} \right]
  &&
 \\ \\
 &=&\log\log n
 &-& \int_1^{\infty} \frac{e^{-y}\rmd y}{y}
 &+& \int_0^1 \frac{1-e^{-y}}{y}\rmd  y
  &+&o(1) \\
 &=&\log\log n &+&& \g & &+&o(1),
 \end{array}
 \end{eqnarray*}
by plugging $z=1$ into~(\ref{eqE1}), taking into account that for
$y\in [0,1]$
 \bn
  1-e^{-B_n(y)}=(1-e^{-y})+\frac{\g ye^{-y}}{\log n}+y\times
  O((\log n)^{-2}),
 \en
and at the same time
 \bn
 \left| \int_1^{\log n} \frac{e^{-B_n(y)}\rmd y}{y} - \int_1^{\infty} \frac{e^{-y}\rmd  y}{y}\right|
 &\le&\int_1^{\sqrt{\log n}} \frac{|e^{-B_n(y)}-e^{-y}|\rmd y}{y}
  \\ \\
    +\int_{\sqrt{\log n}}^{\log n}\frac{e^{-B_n(y)}\rmd y}{y}
    +\int_{\sqrt{\log n}}^{\infty} \frac{e^{-y}\rmd y}{y}
  &=:& (I)+(II)+(III),
 \en
where
 \bn
 (I)&\le &\int_1^{\sqrt{\log n}} \frac{\mbox{Const}}{\sqrt{\log n}}\frac{e^{-y}\rmd y}{y}\le \frac{\mbox{Const}}{\sqrt{\log n}},
 \\ \\
 (II)&=&\int_{(\log n)^{-\frac 12}}^1 \frac{\rmd x}{x(1+x/1)\dots(1+x/n)}
  \le \int_{(\log n)^{-\frac 12}}^1 \frac{\rmd x}{x\left(1+x\left(1+\frac 12+\dots+\frac 1n\right)\right)}
  \\  \\
  &=& \left. \log \frac{x}{1+H_{n}} \right|_{(\log n)^{-\frac 12}}^1
  =\log\frac{H_n+\sqrt{\log n}}{H_n+1} = \frac{1}{\sqrt{\log n}} \, (1+o(1)),\\ \\
 (III)&\le&\frac 1{\sqrt{\log n}}\int_1^{\infty} e^{-y} \rmd y=\frac 1{\sqrt{\log
 n}}.
 \en
\Cox

{\sl Proof of Lemma~\ref{lem_com2}.} The derivation of (a) is
fairly straightforward by induction; it can also be recovered from
Section 4 in~\cite{BRAN}.

To establish (b), note that
$$
a(n,m)\le \sum_{i_1=1}^n  \sum_{i_2=1}^n\dots  \sum_{i_m=1}^n
\frac 1{i_1 i_2\dots i_m}= \left(\sum_{i=1}^n \frac 1i \right)^m
\equiv (H_n)^m
$$
and $H_n\le 1+\log n$ for $n\ge 1$.

Finally, to prove (c), let
 \bn
 \tilde a(n,m) :=\sum_{1\le i_1< i_2<\dots< i_m\le n} \frac 1{i_1 i_2\dots i_m}
 \en
(observe that here all $i_k$'s must be distinct). From
equation~(3.2) in~\cite{GRUN} it follows that for a fixed $m$
satisfy
 \bn
 \tilde a(n,m)=\frac{(\log n)^m}{m!}+\frac{\g (\log n)^{m-1}}{{(m-1)}!}+
 \frac{(\g^2-\zeta(2))\log^{m-2} n}{(m-2)!\, 2}+\dots
 \en
On the other hand,
 \bn
 0&<&a(n,m)-\tilde a(n,m)=\sum_{r=1}^{m-1} \sum_{1\le i_1< i_2<\dots<i_r=i_{r+1}\le i_{r+2}\le\dots\le i_m\le n} \frac 1{i_1 i_2\dots i_m}
\\ \\
&\le& \sum_{r=1}^{m-1}  \left(\sum_{k=1}^n \frac
1{k^2}\right)\sum_{1\le i_1< \dots<i_{r-1}\le i_{r+2}\le\dots\le
i_m\le n} \frac 1{i_1 i_2\dots i_m} \\ \\
 &\le & 2 m\, a(n,m-2)\le 2m (\log n+1)^{m-2}
 \en
Therefore,
$$
a(n,m)=\frac{(\log n)^m}{m!}+\frac{\g (\log
n)^{m-1}}{{(m-1)}!}+O(\log^{m-2} n)
$$
similar to $\tilde a(n,m)$. \Cox

\section*{Acknowledgment}
The author wishes to thank the anonymous referee for helpful
suggestion and corrections, and B\'alint T\'oth for useful
discussions.

\end{document}